\documentclass{amsart}
\usepackage[colorlinks]{hyperref}
\AtBeginDocument{
  \hypersetup{
    linkcolor=blue,
    citecolor=green,
  }
}
\usepackage{cleveref}
\crefname{equation}{}{}

\newtheorem{theorem}{Theorem}[section]
\newtheorem{lemma}[theorem]{Lemma}

\theoremstyle{definition}

\theoremstyle{remark}
\newtheorem{remark}[theorem]{Remark}

\numberwithin{equation}{section}

\begin{document}

\title{A note on the Gagliardo-Nirenberg inequality in a bounded domain}

\author{Congming Li}
\address{School of Mathematical Sciences, Shanghai Jiao Tong University, Shanghai, China}
\email{congming.li@sjtu.edu.cn}
%
\thanks{This research is supported by the National Natural Science Foundation of China (Grant No. 12031012 and 11831003) and the Institute of Modern Analysis-A Frontier Research Center of Shanghai.}

\author{Kai Zhang}
\address{School of Mathematical Sciences, Shanghai Jiao Tong University, Shanghai, China}
\email{zhangkaizfz@gmail.com}
\thanks{}

\subjclass[2020]{Primary 26D10, 35A23 , 46E35}

\date{October 1, 2021}


\keywords{Gagliardo-Nirenberg inequality, Interpolation inequality, bounded domain}

\begin{abstract}
The classical Gagliardo-Nirenberg inequality was established in $\mathbb{R}^n$. An extension to a bounded domain was given by Gagliardo in 1959. In this note, we present a simple proof of this result and prove a new Gagliardo-Nirenberg inequality in a bounded Lipschitz domain.
\end{abstract}

\maketitle

\section{Introduction}\label{S1}

In this short note, we prove two kinds of the Gagliardo-Nirenberg inequalities in a bounded domain based on the Gagliardo-Nirenberg inequality in $\mathbb{R}^n$. First, we introduce some notations. Denote by $|E|$ the Lebesgue measure of $E\subset \mathbb{R}^n$. Let $\Omega \subset \mathbb{R}^n$ be a domain and $u:\Omega\rightarrow \mathbb{R}$. For $0<\alpha\leq 1$, define
\begin{equation*}
[u]_{C^{\alpha}(\bar\Omega)}=\sup_{x,y\in \Omega,x\neq y} \frac{|u(x)-u(y)|}{|x-y|^{\alpha}}.
\end{equation*}
For $0<p\leq \infty$, denote
\begin{equation*}
\|u\|_{p,\Omega}=\left(\int_{\Omega} |u|^p\right)^{1/p}.
\end{equation*}
For $-\infty<p\leq +\infty$, set
\begin{equation*}
|u|_{p,\Omega}=\left\{\begin{aligned}
&\|u\|_{p,\Omega} &&p>0;\\
&\|\nabla^k u\|_{\infty,\Omega} &&p<0, -n/p=k\in \mathbb{N};\\
&[\nabla^k u]_{C^{\alpha}(\bar\Omega)} &&p<0, -n/p=k+\alpha, k\in \mathbb{N}, 0<\alpha<1.
  \end{aligned}\right.
\end{equation*}
If $|u|_{p,\Omega}<\infty$, we say that $u\in L^{p}(\Omega)$. For simplicity, we also write $\|u\|_{p}$ and $|u|_{p}$ instead of $\|u\|_{p,\Omega}$ and $|u|_{p,\Omega}$ if $\Omega$ is clearly understood.

The famous Gagliardo-Nirenberg inequality was first proved by Gagliardo \cite{MR109295} and Nirenberg \cite{MR0109940} independently restricted in Sobolev spaces $W^{k,p}(\mathbb{R}^n)$ where $k\in \mathbb{N}$ (i.e. nonnegative integers) and $1\leq p\leq \infty$. More precisely (see \cite[Theorem 12.87]{MR3726909}),
\begin{theorem}\label{th1.0}
Let $1\leq q,r\leq +\infty$, $k,j\in \mathbb {N}$ with $j<k$, $j/k\leq \theta \leq 1$ and $p\in \mathbb{R}$ such that
\begin{equation}\label{e1.0}
\frac{n}{p}-j=\theta \left(\frac{n}{r}-k\right)+(1-\theta)\frac{n}{q}.
\end{equation}
Then there exists a constant $C$ depending only on $n,k,q,r$ and $\theta$ such that for any $u\in W^{k,r}(\mathbb{R}^n)\cap L^{q}(\mathbb{R}^n)$,
\begin{equation}\label{e1.0-2}
|\nabla^{j} u|_{p}\leq C\|\nabla^{k} u\|_{r}^{\theta}\|u\|_{q}^{1-\theta}
\end{equation}
with the exception that if $1<r<+\infty$ and $k-j-n/r\in \mathbb{N}$, we must take $j/k\leq \theta <1$.
\end{theorem}
\begin{remark}\label{re1.1}
The relation \cref{e1.0} is a necessary requirement of the scaling consideration in \cref{e1.0-2}.
\end{remark}

The Gagliardo-Nirenberg inequality has been extended greatly in different directions, such as inequalities in other function spaces \cite{MR2424108, MR2524231, MR4135294, MR2215922}, involving fractional derivatives \cite{MR3813967,MR3990737}, aim at best constants \cite{MR1940370, MR591555}, in bounded domains\cite{MR4049810,MR3813967,MR3990737} and manifolds \cite{MR2456207,MR2381156} ect.

Our main results are the following:
\begin{theorem}\label{th1.1}
Let $\Omega$ be a bounded Lipschitz domain. Assume that $1\leq q,r\leq +\infty$, $k,j\in \mathbb {N}$ with $j<k$, $j/k\leq \theta \leq 1$ and $p\in \mathbb{R}$ such that
\begin{equation}\label{e1.1}
\frac{n}{p}-j=\theta \left(\frac{n}{r}-k\right)+(1-\theta)\frac{n}{q}.
\end{equation}
Then there exists a constant $C$ depending only on $n,k,q,r,\theta$ and $\Omega$ such that for any $u\in W^{k,r}(\Omega)\cap L^{q}(\Omega)$,
\begin{equation}\label{e1.2}
|\nabla^{j} u|_{p}\leq C\|\nabla^{k} u\|_{r}^{\theta}\|u\|_{q}^{1-\theta}+C\|u\|_{q}
\end{equation}
with the exception that if $1<r<+\infty$ and $k-j-n/r\in \mathbb{N}$, we must take $j/k\leq \theta <1$.
\end{theorem}
\begin{remark}\label{re1.2}
\Cref{th1.1} is proved in \cite{MR109295} and stated in \cite{MR208360}. Here, we give a simple proof.
\end{remark}

\begin{theorem}\label{th1.2}
Let $\Omega$ be a bounded Lipschitz domain. Assume that $1\leq q,r\leq +\infty$, $k,j\in \mathbb {N}$ with $j<k$, $j/k\leq \theta \leq 1$ and $p\in \mathbb{R}$ such that
\begin{equation}\label{e1.5}
\frac{n}{q}>\frac{n}{r}-k
\end{equation}
and
\begin{equation*}\label{e1.3}
\frac{n}{p}-j=\theta \left(\frac{n}{r}-k\right)+(1-\theta)\frac{n}{q}.
\end{equation*}
Then for any $E\subset \Omega$ with $|E|>0$, there exists a constant $C$ depending only on $n,k,q,r,\theta, E$ and $\Omega$ such that for any $u\in W^{k,r}(\Omega)$,
\begin{equation}\label{e1.4}
|\nabla^{j} u-\nabla^{j}P_{E}|_{p}\leq C\|\nabla^{k} u\|_{r}^{\theta}\|u-P_{E}\|_{q}^{1-\theta}
\end{equation}
with the exception that if $1<r<+\infty$ and $k-j-n/r\in \mathbb{N}$, we must take $j/k\leq \theta <1$.
The $P_E$ is the unique polynomial of degree $k-1$ such that for any multi-index $0\leq |\gamma|\leq k-1$,
\begin{equation*}
  \int_{E} \nabla^{\gamma}u=\int_{E} \nabla^{\gamma}P_{E}.
\end{equation*}
\end{theorem}
\begin{remark}\label{re1.3}
The \cref{e1.5} is a technical assumption and we don't know whether it can be removed or not.
\end{remark}

\section{Proof of the main results}
To prove \Cref{th1.1}, we need the following two lemmas. The first lemma (see \cite{MR2225465} and \cite[Theorem 5, Chapter IV]{MR0290095}) deals with extension from a bounded domain to the whole space. This allows us to use the Gagliardo-Nirenberg inequality in the whole space.
\begin{lemma}\label{le2.1}
Let $\Omega\subset \mathbb{R}^n$ be a bounded domain. Then there exists a bounded linear extension operator $T:W^{k,p}(\Omega)\rightarrow W^{k,p}(\mathbb{R}^n)$ ($k\geq 0, 1\leq p\leq \infty$) such that
\begin{equation}\label{e2.1}
\|Tf\|_{W^{k,p}(\mathbb{R}^n)}\leq C\|f\|_{W^{k,p}(\Omega)},~~\forall f\in W^{k,p}(\Omega),
\end{equation}
where $C$ depends only on $n,k,p$ and $\Omega$.
\end{lemma}

The next is an interpolation in a bounded domain.
\begin{lemma}\label{le2.2}
Let $\Omega$ be a bounded Lipschitz domain and $u\in W^{k,r}(\Omega)$ ($k\geq 0,1\leq r\leq \infty$). Then for any $\varepsilon,q>0$, $0\leq j<k$ and $p\in \mathbb{R}$ with
\begin{equation}\label{e2.6}
\frac{n}{p}-j> \frac{n}{r}-k,
\end{equation}
there exists a constant $C$ depending only on $n,k,q,r,\varepsilon$ and $\Omega$ such that
\begin{equation*}
|\nabla^{j} u|_{p}\leq \varepsilon\|\nabla^{k} u\|_{r}+C\|u\|_{q}.
\end{equation*}
\end{lemma}
\proof We prove the lemma by contradiction. Suppose that the conclusion is false. Then there exist $\varepsilon_0>0$ and a sequence of $\left\{u_m\right\}\subset W^{k,r}(\Omega)$ such that
\begin{equation*}
|\nabla^{j} u_m|_{p}\geq \varepsilon_0\|\nabla^{k} u_m\|_{r}+m\|u_m\|_{q}.
\end{equation*}

Without loss of generality, we assume that $\|u_m\|_{W^{k-1,r}(\Omega)}+|\nabla^j u_m|_{p}=1$. Then
\begin{equation}\label{e2.2}
  1\geq \varepsilon_0\|\nabla^{k} u_m\|_{r}+m\|u_m\|_{q}.
\end{equation}
Hence, $\|u_m\|_{W^{k,r}(\Omega)}\leq 1+1/\varepsilon_0$. By \cref{e2.6}, $W^{k,r}$ is compactly embedded into $L^{p}$. Thus, there exist a subsequence (denoted by $\left\{u_m\right\}$ again) and $\bar{u}\in W^{k,r}(\Omega)$ such that as $m\rightarrow \infty$,
\begin{equation*}
\|u_m-\bar u\|_{W^{k-1,r}(\Omega)}+|\nabla^j u_m-\nabla^j \bar{u}|_{p}\rightarrow 0.
\end{equation*}
Therefore, $\|\bar u\|_{W^{k-1,r}(\Omega)}+|\nabla^j \bar u|_{p}=1$. However, from \cref{e2.2},
\begin{equation*}
  \|\bar u\|_{q}=\lim_{m\rightarrow \infty} \|u_m\|_{q}=0.
\end{equation*}
That is, $\bar{u}\equiv 0$ and we arrive at a contradiction. \qed~\\

Now, we can give the

\noindent\textbf{Proof of \Cref{th1.1}.} Extend $u$ to $Tu\in W^{k,r}(\mathbb{R}^n)\cap L^q(\mathbb{R}^n)$. By the Gagliardo-Nirenberg inequality in $\mathbb{R}^n$ (see \Cref{th1.0}) and \Cref{le2.2} with $\varepsilon=1$,
\begin{equation*}
  \begin{aligned}
    |\nabla^j u|_{p,\Omega}&\leq |\nabla^j Tu|_{p,\mathbb{R}^n}\\
    &\leq C\|\nabla^k Tu\|^{\theta}_{r,\mathbb{R}^n}\|Tu\|^{1-\theta}_{q,\mathbb{R}^n}\\
    &\leq C\left(\sum_{i=1}^{k}\|\nabla^iu\|_{r,\Omega}\right)^{\theta}
    \|u\|^{1-\theta}_{q,\Omega}\\
    &\leq C(\|u\|_{q,\Omega}+\|\nabla^ku\|_{r,\Omega})^{\theta}
    \|u\|^{1-\theta}_{q,\Omega}\\
    &\leq C\|\nabla^ku\|^{\theta}_{r,\Omega}\|u\|^{1-\theta}_{q,\Omega}+
    C\|u\|_{q,\Omega}.
  \end{aligned}
\end{equation*}
\qed~\\

To prove \Cref{th1.2}, we need the following lemma.
\begin{lemma}\label{le2.3}
Let $\Omega$ be a bounded Lipschitz domain, $E\subset \Omega$ with $|E|>0$ and $u\in W^{k,r}(\Omega)$ ($k\geq 0,1\leq r\leq \infty$). Then there exists a unique polynomial $P_{E}$ of degree $k-1$ with
\begin{equation}\label{e2.3}
\int_{E} \nabla^{\gamma}u=\int_{E} \nabla^{\gamma} P_E,~~\forall~~ 0\leq |\gamma|\leq k-1
\end{equation}
such that for any $0\leq j<k$ and $p\in \mathbb{R}$ with
\begin{equation*}
\frac{n}{p}-j> \frac{n}{r}-k,
\end{equation*}
we have
\begin{equation*}
|\nabla^{j} u-\nabla^{j}P_E|_{p}\leq C\|\nabla^{k} u\|_{r},
\end{equation*}
where $C$ depends only on $n,k,r,\Omega$ and $E$.
\end{lemma}
\begin{remark}\label{re2.1}
\Cref{le2.3} can be regarded as an extension of the Poincar\'{e} inequality.
\end{remark}
\proof The proof is similar to that of \Cref{le2.2}. Suppose that the conclusion is false. Then there exist a sequence of $\left\{u_m\right\}\subset W^{k,r}(\Omega)$ and $\left\{P_m\right\}$ such that \cref{e2.3} holds and
\begin{equation*}
|\nabla^{j} u_m-\nabla^{j} P_m|_{p}\geq m\|\nabla^{k} u_m\|_{r}.
\end{equation*}

Let $v_m=u_m-P_m$. Without loss of generality, we assume that $\|v_m\|_{W^{k-1,r}(\Omega)}+|\nabla^j v_m|_{p}=1$. Then
\begin{equation}\label{e2.4}
  1\geq m\|\nabla^{k} v_m\|_{r}.
\end{equation}
Hence, $\|v_m\|_{W^{k,r}(\Omega)}\leq 2$. By the compact embedding, there exist a subsequence (denoted by $\left\{v_m\right\}$ again) and $\bar{v}\in W^{k,r}(\Omega)$ such that as $m\rightarrow \infty$,
\begin{equation*}
\|v_m-\bar v\|_{W^{k-1,r}(\Omega)}+|\nabla^j v_m-\nabla^j \bar{v}|_{p}\rightarrow 0.
\end{equation*}
Thus,
\begin{equation}\label{e2.5}
\|\bar v\|_{W^{k-1,r}(\Omega)}+|\nabla^j \bar v|_{p}=1.
\end{equation}
However, from \cref{e2.4},
\begin{equation*}
  \|\nabla^k\bar v\|_{r}=0.
\end{equation*}
That is, $\bar v$ is a polynomial of degree $k-1$. From \cref{e2.3}, $\int_E \nabla^{\gamma} v_m=0$ for any $m\geq 1$ and $0\leq |\gamma|\leq k-1$. Let $m\rightarrow \infty$, we have the same equalities for $\bar{v}$. Hence $\bar{v}\equiv 0$, which contradicts with \cref{e2.5}. \qed~\\

Now, we can give the

\noindent\textbf{Proof of \Cref{th1.2}.} Since $n/q>n/r-k$, by \Cref{le2.3},
\begin{equation*}
\|u-P_{E}\|_{q}\leq C \|\nabla^{k} u\|_{r}.
\end{equation*}
Hence, from \Cref{th1.1},
\begin{equation*}
|\nabla^{j} u-\nabla^{j}P_{E}|_{p}\leq C\|\nabla^{k} u\|_{r}^{\theta}\|u-P_{E}\|_{q}^{1-\theta}+\|u-P_{E}\|_{q}\leq C\|\nabla^{k} u\|_{r}^{\theta}\|u-P_{E}\|_{q}^{1-\theta}.
\end{equation*}
\qed~\\

\end{document}